\documentclass[reqno,11pt]{amsart}
\usepackage{amsmath, latexsym, amsfonts, amssymb, amsthm, amscd}
\usepackage{epsfig}
\usepackage{graphics,epsf,psfrag}

\numberwithin{equation}{section}

\newtheorem{theorem}{Theorem}[section]
\newtheorem{lemma}[theorem]{Lemma}
\newtheorem{prop}[theorem]{Proposition}

\newtheorem{rem}[theorem]{Remark}

\newcommand{\R}{\mathbb{R}}
\newcommand{\Z}{\mathbb{Z}}
\newcommand{\N}{\mathbb{N}}
\renewcommand{\tilde}{\widetilde}

\newcommand{\cA}{{\ensuremath{\mathcal A}} }
\newcommand{\cF}{{\ensuremath{\mathcal F}} }

\newcommand{\cC}{{\ensuremath{\mathcal C}} }
\newcommand{\cN}{{\ensuremath{\mathcal N}} }
\newcommand{\cL}{{\ensuremath{\mathcal L}} }

\newcommand{\cV}{{\ensuremath{\mathcal V}} }

\newcommand{\bP}{{\ensuremath{\mathbf P}} }
\newcommand{\bE}{{\ensuremath{\mathbf E}} }

\newcommand{\bbE}{{\ensuremath{\mathbb E}} }

\newcommand{\bbN}{{\ensuremath{\mathbb N}} }

\newcommand{\bbP}{{\ensuremath{\mathbb P}} }

\newcommand{\bbR}{{\ensuremath{\mathbb R}} }

\renewcommand{\tilde}{\widetilde}          
\DeclareMathSymbol{\leqslant}{\mathalpha}{AMSa}{"36} 
\DeclareMathSymbol{\geqslant}{\mathalpha}{AMSa}{"3E} 
\DeclareMathSymbol{\eset}{\mathalpha}{AMSb}{"3F}     
\newcommand{\dd}{\text{\rm d}}             

\newcommand{\ga}{\alpha}

\newcommand{\gd}{\delta}
\newcommand{\gep}{\varepsilon}       

\newcommand{\gs}{\sigma}

\newcommand{\p}{\bbP}

\def\bs{\boldsymbol}
\def\ind{\bs{1}}

\newcommand{\rc}{\mathrm c}
\newcommand{\rf}{\mathrm f}

\newcommand{\tf}{\textsc{f}}

\title[Wetting models in (1+1)--dimension]{Sharp asymptotic behavior\\ for wetting models in (1+1)--dimension}

\author{Francesco Caravenna}
\address{Institut f\"ur Mathematik, Universit\"at Z\"urich,
Winterthurerstrasse 190, CH-8057 Z\"urich}
\email{francesco.caravenna\@@math.unizh.ch}
\urladdr{http://www.matapp.unimib.it/\raisebox{0.11ex}{\tiny$\sim$}fcaraven/}

\author{Giambattista Giacomin}
\address{Laboratoire de Probabilit{\'e}s de P 6\ \& 7 (CNRS U.M.R. 7599) and
Universit{\'e} Paris 7 -- Denis Diderot, U.F.R. Mathematiques, Case 7012, 2 place Jussieu, 75251 Paris cedex 05, France}
\email{giacomin\@@math.jussieu.fr}
\urladdr{http://www.proba.jussieu.fr/pageperso/giacomin/GBpage.html}

\author{Lorenzo Zambotti}
\address{Dipartimento di Matematica, Politecnico di Milano,
Piazza Leonardo da Vinci 32, 20133 Milano, Italy}
\email{lorenzo.zambotti\@@polimi.it}
\urladdr{http://www1.mate.polimi.it/\raisebox{0.11ex}{\tiny$\sim$}zambotti/}

\date{\today}



\begin{document}

\begin{abstract}
We consider continuous and discrete (1+1)--dimensional wetting models
which undergo a localization/delocalization
phase transition. Using a simple approach based on Renewal Theory
we determine the precise asymptotic behavior of the partition function,
from which we obtain the scaling limits of the models and
an explicit construction of the infinite
volume measure (thermodynamic limit) in all regimes, including
the critical one.
\end{abstract}

\maketitle

\section{Introduction}

\subsection{Definition of the model}

The building blocks of our model are a $\gs$--finite measure $\mu$ on $\R$
(the {\sl single site a priori measure}) and a function $V:\R \mapsto \R\cup \{+\infty\}$
(the {\sl potential}). We allow two possible choices of~$\mu$:
\begin{itemize}
\item {\bf Continuous set--up:} $\mu= \dd x$ is the Lebesgue
measure on~$\R$. In this case we require that $\exp (-V(\cdot))$
be continuous.
\item \rule{0pt}{1.2em}{\bf Discrete set--up:} $\mu$ is the counting measure on $\Z\,$.
\end{itemize}
In both settings we assume that $V(0) < \infty$ and that
\begin{equation*}
\kappa \, := \, \int_{\mathbb R} e^{-V(y)} \,\mu( \dd y) \, <\, \infty \,.
\end{equation*}
Additional assumptions on $V(\cdot)$ will be stated in the next subsection.

\begin{figure}[h]
\begin{center}
\leavevmode
\epsfysize=5 cm
\psfragscanon
\psfrag{S}[c][c]{\Large $S_n$}
\psfrag{n}[c][c]{ \Large $n$}
\psfrag{0}[c][c]{$0$}
\psfrag{N+1}[c][c]{ $N+1$}
\epsfbox{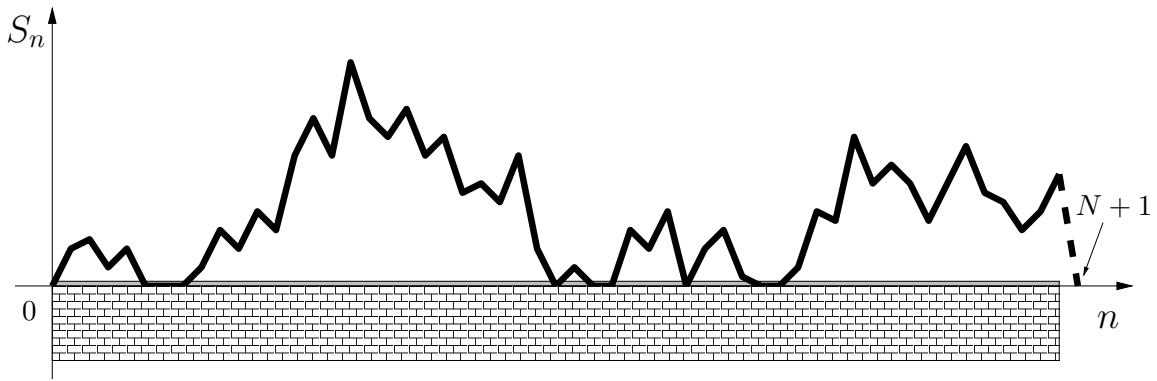}
\end{center}
\caption{\label{fig:1}
A trajectory of the walk (or interface) above the wall. The last step
may be constrained or free, according to the model.
The walk is rewarded when it enters the gray {\sl thin} layer close to the wall. The competition between the {\sl energetic gain} coming from
this reward and the {\sl entropic repulsion} due to the
presence of the wall leads to a nontrivial behavior.
Such a model has been proposed at several instances
for the study of interfaces and polymers, we refer
to \cite{cf:Yvan} for details and references.
}
\end{figure}

\smallskip
For $\gep \ge 0,\ N \in \N$ our model is defined by the following probability measure
on $(\R^+)^N := [0,\infty)^N\,$:
\begin{equation*}
{\bf P}^a_{\gep ,N}(dx) \, := \, \frac{1}{Z^a_{\gep ,N}} \,
\exp(-H_N^a(x)) \, \prod_{i=1}^N
\big(1_{(x_i>0)} \, \mu (dx_i) +\gep \delta_0(dx_i) \big) \,,
\end{equation*}
where $Z_{\gep, N}^a$ is the normalizing constant ({\sl partition function}),
$a$ is a label that stands for $\rf$ ({\sl free}) or $\rc$
({\sl constrained}) and the corresponding Hamiltonians are defined by
\begin{gather*}
H_N^\rf(x) \, := \, \sum_{i=0}^{N-1} V(x_{i+1}-x_i)\, \qquad x_0 :=0 \,, \\
H_N^\rc(x) \, := \, \sum_{i=0}^{N} V(x_{i+1}-x_i)\, \qquad x_0 := x_{N+1} := 0\,.
\end{gather*}
We interpret~${\bf P}^a_{\gep ,N}$ as an effective model for
a $(1+1)$--dimensional interface above an impenetrable wall that,
when~$\gep > 0$, attracts it (see Figure~\ref{fig:1} ad the relative caption).


\subsection{A random walk viewpoint}
We introduce a sequence $(Y_i)_{i\in\N}$ of IID random variables with law $\bP$ such that
\begin{equation*}
    \bP (Y_1 \in \dd x) \;:=\; \frac1 {\kappa} \;
    {\exp(-V(x))} \; \mu(\dd x)\,,
\end{equation*}
and we denote by $(S_n)_{n\ge 0}$ the associated random walk: $S_0:=0$,
$S_n:=Y_1+\ldots+Y_n$.

The basic assumption we make on the potential~$V(\cdot)$ is the following one:
\begin{itemize}
\item[(H)] \rule{0pt}{1.1em}The {\sl truncated variance} $t\mapsto \cV(t) := \bE \big[ |Y_1|^2 \ind_{\{|Y_1| \le t\}} \big]$ is slowly varying at infinity and $\bE\big[Y_1\big] = 0\,$.
\end{itemize}
We recall that a function $L(\cdot)$ is said to be slowly varying at infinity if
for every $c>0$ one has $L(ct)/L(t) \to 1$ as $t\to\infty$. This entails that
$L(x)/x^{\ga} \to 0$ as $x\to\infty$, for every $\ga > 0$, cf. \cite[Prop.~1.3.6]{cf:BinGolTeu}.
Notice that if the truncated variance is slowly varying, then we have
$\bP \big[ |Y_1| \ge t \big] \le \cV(t)/t^2 $ for large~$t$, cf. \cite[Th.~8.3.1]{cf:BinGolTeu},
hence all the moments of $Y_1$ of order less than~$2$ are automatically finite.

Of course assumption~(H) holds whenever $Y_1$ is centered and has a {\sl finite variance},
that is when $\bE\big[ Y_1 \big] = 0$ and $\bE \big[ |Y_1|^2 \big] =:\gs^2 < \infty\,$.
Allowing the truncated variance to be slowly varying turns out to be a
very natural generalization: indeed assumption~(H) is a {\sl necessary and sufficient}
condition for $(S_n)_{n\ge 0}$ to be in the domain of attraction (without centering)
of the Gaussian law, see Appendix~\ref{app:1} for more on this issue (cf. also
\cite[Th.~8.3.1]{cf:BinGolTeu}).

Now let us look more closely to our model. For $\gep = 0$ we have
the following random walk interpretation: $\bP_{0,N}^\rf$ is just the law of $( S_1, \ldots, S_N )$
under the positivity constraint $\{S_1> 0, S_2> 0, \ldots, S_N > 0\}$, while
$\bP_{0,N}^\rc$ is the law of the same random vector under
the further constraint $\{S_{N+1}=0\}$. Then by the weak convergence
towards Brownian meander and Brownian excursion we have that the
$\bP_{0,N}^a$--typical height of the interface in the bulk is of
order~$\sqrt{N}$, hence {\sl very far} from the interface  (\textsl{delocalized regime}).
On the other hand, when $\gep>0$ the interface receives an $\gep$--reward
each time it touches the wall and intuitively
one expects that if $\gep$ is large enough this
attractive effect should be able to beat the entropic repulsion, leading to a
\textsl{localized regime}. As we are going to see, this scenario is
correct.

\begin{rem}\rm
We stress that the assumption $V(0) < \infty$ has been made only for simplicity and
can be removed at the price of some heavier notation. Also the condition that
$\exp (-V(\cdot))$ be a continuous function (in the continuous set--up) can be lightened
to the requirement that for some~$n \in \N$ the $n$ fold convolution
of $\exp (-V(\cdot))$ with itself be continuous, in order to apply the Local Limit Theorem
for Densities (see Appendix~\ref{app:1}).
\end{rem}

\subsection{The phase diagram and the scaling limits}

The model we are considering has been studied in~\cite{cf:IY} in the discrete set--up,
for the special choice $\exp(-V(x)) = \frac v2 \, \ind_{(|x| = 1)} + r\, \ind_{(x=0)}$, with $v>0$, $r\ge 0$ and $v+r=1$,
and more recently in~\cite{dgz} in the finite variance continuous
set--up, that is when
$\bE\big[ |Y_1|^2 \big] < \infty$. In both cases it has been proven that:
\begin{itemize}
\item[(i)] There is a phase transition at $\gep = \gep_c > 0$,
between a {\sl delocalized regime} ($\gep \le \gep_c$) in which
the interface is repelled by the wall and a {\sl localized regime} ($\gep > \gep_c$)
in which the interface sticks close to the wall. A convenient definition of (de)localization
may be given for instance in terms of the free energy, that is by looking at the Laplace
asymptotic behavior of the partition function, cf.~\cite[\S~2.2]{dgz}.
\item[(ii)] \rule{0pt}{13pt}More quantitatively, Brownian scaling limits hold, inducing
a further distinction in the delocalized regime. More precisely,
the linearly interpolated diffusive rescaling of $\bP_{\gep,N}^a$ converges in distribution
as $N \to \infty$:
\begin{itemize}
\item \rule{0pt}{12pt}when $\gep < \gep_c$ ({\sl strictly delocalized regime}),
to the Brownian meander if $a = \rf$ or to
the normalized Brownian excursion if $a = \rc$;
\item \rule{0pt}{12pt}when $\gep = \gep_c$ ({\sl critical regime}),
to the reflecting Brownian motion if $a = \rf$
or to the reflecting Brownian bridge if $a = \rc$;
\item \rule{0pt}{12pt}when $\gep > \gep_c$ ({\sl localized regime}),
to the law concentrated on the function taking
the constant value~$0$ for both $a = \rf$ and $a = \rc$.
\end{itemize}
\end{itemize}
The proof of these results in~\cite{cf:IY} has been obtained by exploiting
some very peculiar properties enjoyed by walks with increments in $\{-1, 0,+1\}$. On the other hand,
the more general approach adopted in~\cite{dgz} is based on bounds on the asymptotic
behavior of the partition function~$Z^a_{\gep,N}$ as~$N\to\infty$.

\subsection{Outline of the results}

The purpose of this note is to present a simple approach
based on Renewal Theory, which is applicable in complete generality in both
the continuous and discrete set--ups, that allows to determine the
\textsl{precise asymptotic behavior} of the partition function~$Z^a_{\gep,N}$
in all regimes. This moreover yields a considerable simplification of
several steps in~\cite{dgz}, allowing the extension
of the above results (i) and~(ii) to the general continuous and discrete set--ups
in a straightforward way.

Another important byproduct of our approach,
and possibly the main result
presented here,    concerns the weak convergence
as $N\to\infty$ of $\bP_{\gep,N}^a$ as a probability measure
on $(\R^+)^\N$, the so--called {\sl thermodynamic limit}. This issue
has been already considered in~\cite{dgz}, in the finite variance
continuous set--up and only for the localized regime.
Here we show that the weak limit of~$\bP_{\gep,N}^a$
exists, in both the continuous and discrete set--ups, for  $a=\rf$ and
$a = \rc$ and for all values of~$\gep$, cf. Theorem~\ref{th:main}.
This will come with an explicit description of the limit measure,
whose properties differ considerably in the strictly
delocalized, critical and localized regimes,
in complete analogy to the above mentioned scaling limits.

The exposition is organized as follows:
\begin{itemize}
\item[--] In Section~\ref{sec:ren_th} we
describe a Renewal Theory approach to our model.
\item[--] \rule{0pt}{1.1em}This will lead
 to the determination of the {\sl precise} asymptotic behavior
of the partition function (Section~\ref{sec:sharp_as}), to be compared to \cite[Lemma 3]{dgz}.
\item[--] \rule{0pt}{1.1em}These results, in turn, are the key to proving the
existence of the thermodynamic limit in Section~\ref{sec:inf_vol}
in all regimes.
\item[--] \rule{0pt}{1.1em}Finally,
in Section~\ref{sec:scaling_lim} we give the main ingredients to extend
the proof of the scaling limits given in~\cite{dgz} to our general setting.
\end{itemize}

\section{A Renewal Theory viewpoint}
\label{sec:ren_th}

In this section we make explicit the link with Renewal
Theory, showing that a suitable modification of the constrained partition
function $Z_{\gep,N}^\rc$ can be interpreted as the (generalized) Green function
associated with a probability measure $q(\cdot)$ that we define below.
We also give a simple relation linking $Z_{\gep,N}^\rf$ and $Z_{\gep,N}^\rc$.
Throughout the paper, we write for positive sequences $(a_n)_n$ and $(b_n)_n$:
\[
a_n \sim b_n \ \Longleftrightarrow \  \lim_{n\to\infty} a_n / b_n =1.
\]

\subsection{The constrained case}

Let us consider first the $\gep = 0$ case. We claim that
\begin{equation}\label{eqq}
\frac{Z^\rc_{0 ,N}}{\kappa^{N+1}} \, \sim \, \frac C{\sqrt{2\pi}} \, \frac{L(N)}{N^{3/2}}
 \qquad (N\to\infty),
\end{equation}
where $C$ is a positive constant and $L(\cdot)$ is a slowly varying function
given in Appendix~\ref{app:1}.
The proof of this relation is deferred to Appendix~\ref{app:asympt}. We stress that
in the case of finite variance we have $L(\cdot) \equiv \gs^{-1}$.

We set for convenience $Z^\rc_{\gep,0} := \exp(-V(0))$. We recall that $L(x)/x^\ga \to 0$
as $x\to\infty$ for every $\ga > 0$, because $L(\cdot)$ is slowly varying.
Then from equation \eqref{eqq} it follows that
$\gamma :=\sum_{n\ge 1} Z^\rc_{0 ,n-1}/\kappa^{n}<\infty$, hence we can define a probability
distribution $q(\cdot)$ on~$\N$ by setting
\[
q(n) \, := \, \frac 1 \gamma \, \frac{Z^\rc_{0 ,n-1}}{\kappa^n} \qquad (n\geq 1)\,,
\]
and we have for some positive constant~$c_q$:
\begin{equation}
\label{c_q}
q(n) \, \sim \, c_q \, \frac{L(n)}{n^{3/2}} \qquad \quad (n\to\infty)\,.
\end{equation}

Next we pass to the $\gep>0$ case. It is convenient to switch from the parameter $\gep$
to $\gd := \gamma \gep$ and to make a change of scale, by setting for~$N\ge 1$
\begin{equation}\label{eq:zeta_tilde}
\widetilde Z^\rc_{\delta,N} \, := \, \frac{(\gd/\gamma) \, Z^\rc_{(\gd/\gamma)
,N-1}}{\kappa^{N}} \qquad \quad
\widetilde Z^\rf_{\delta,N} \, := \, \frac{(\gd/\gamma) \, Z^\rf_{(\gd/\gamma)
,N}}{\kappa^{N}}\,.
\end{equation}
Then from the very definition of the partition function~$Z_{\gep,N}^c$
it follows that $\widetilde Z^\rc_{\delta,N+1}$ satisfies
the following recurrence relation (cf.~\cite[Lemma~2]{dgz}):
\begin{equation}\label{e2}
\widetilde Z^\rc_{\delta,0} \, := \, 1\,, \qquad \ \
\widetilde Z^\rc_{\delta,N} \, = \, \delta\sum_{t=1}^{N} q(t)
\, \widetilde Z^\rc_{\delta,N-t}\,, \quad (N\geq 1)\,.
\end{equation}
This is nothing but the \textsl{(generalized) renewal equation} driven
by~$\delta \, q(\cdot)$, see~\cite{asmu}, and it is easily checked
that its unique solution is given by
\begin{equation} \label{eq:ren_sol}
    \widetilde Z^\rc_{\delta,N} \, = \, \sum_{k=0}^\infty \delta^k \, q^{k*}(N),
\end{equation}
where $q^{k*}$ denotes the $k$--fold convolution of $q$ with itself
(by convention $q^{0*}(n):=\ind_{(n=0)}$).
Notice that the infinite sum in the right hand side of
\eqref{eq:ren_sol} is in fact a sum from $0$ to~$N$.

\smallskip
We point out that this approach can be generalized via the so--called
Markov Renewal Theory~\cite{asmu}, allowing to study
{\sl periodically inhomogeneous} models, that is the case in which $\gep$
is substituted by $\gep_i\in \R$ and $\gep_i=\gep_{i+T}$
for some $T\in \N$ and every $i$.
This has been recently worked out in~\cite{cgz},
in the context of models of {\sl copolymers with adsorption}.

\subsection{The free case.}

From the definition of~$ Z_{\gep, N}^\rf$ we have the following simple relation
for the modified free partition function:
\begin{equation}
\label{rf}
\tilde Z_{\gd, N}^\rf \, = \, \sum_{t=0}^{N}
\tilde Z_{\gd, t}^\rc \,  P(N-t)\,,
\end{equation}
where $P(n)$ is the probability that the unperturbed random walk $(S_i)_i$ stays
positive up to epoch~$n$:
\begin{equation}
\label{P}
P(n) \, := \, \bP(S_i> 0, \, i=1,\ldots,n)\,,
\qquad P(0) \, := \, 1\,.
\end{equation}
It is worth recalling that if we set
\begin{equation}
\label{c_Y}
L'(n) \, :=\,
\sqrt{n} \, P (n)\,,
\end{equation}
then assumption~(H) yields that $L'(\cdot)$ is slowly varying at infinity,
cf. for instance~\cite{cf:DonGre}.
We note that by the standard theory of
stability~\cite{cf:BinGolTeu} this is equivalent to saying that
the first descending ladder epoch~$\overline T_1$ (see Appendix~\ref{app:discrete})
is in the domain of attraction of the positive stable law of index~$1/2$, like in
the simple random walk case.

\section{Sharp asymptotic behavior of the partition function}
\label{sec:sharp_as}

In this section we specialize the theory of the
renewal equation~\cite{asmu} to our {\sl heavy--tailed setting}, see~\eqref{c_q},
in order to find the asymptotic behavior of
$\widetilde Z_{\gd, N}^\rc$ and~$\widetilde Z_{\gd, N}^\rf$.

\subsection{The strictly delocalized regime.}
The following lemma gives the asymptotic behavior of the partition function
in the strictly delocalized regime. To get some intuition we observe that
$\widetilde Z^\rc_{\gd,N}$ when $\delta<1$ is the Green function of a renewal
process with \textsl{defective} interarrival distribution $\delta \, q(\cdot)$,
as it follows from~\eqref{eq:ren_sol}.
\begin{lemma}\label{lemma2}
If $\delta<1$ and relations \eqref{c_q} and \eqref{c_Y} hold, then
\begin{gather}\label{e3}
\widetilde Z^\rc_{\delta,N} \, \sim \, \frac{\delta \, c_q}{(1-\delta)^2}
\; \frac{L(N)}{N^{3/2}} \qquad \quad (N\to\infty)\,,\\
\label{e32}
\tilde Z_{\delta, N}^\rf \, \sim \, \frac{1}{1-\delta} \; \frac{L'(N)}{N^{1/2}}
\qquad \quad (N\to\infty)\,.
\end{gather}
\end{lemma}
\noindent
{\bf Proof.}
From~\eqref{eq:ren_sol} we have that
\begin{equation} \label{eq:rel1}
\frac{n^{3/2}}{L(n)} \,\widetilde Z^\rc_{\delta,n} \, = \,
\sum_{k=1}^\infty \delta^k\, \frac{n^{3/2}}{L(n)}\, q^{k*}(n) \,.
\end{equation}
We claim that
\begin{equation}\label{eq1}
\lim_{n\to\infty}\, \frac{n^{3/2}}{L(n)}\, \, q^{k*}(n) \, = \, k \, c_q, \qquad \forall \
k\geq 1.
\end{equation}
We argue by induction, the case $k=1$ being true by
\eqref{c_q}.  Suppose that we have proven \eqref{eq1} for
$k=1, \ldots, m$, then we have:
\begin{eqnarray*}
& &
\frac{n^{3/2}}{L(n)} \, q^{(m+1)*}(n) \, = \, \left[ \sum_{i=1}^{\lfloor n/2 \rfloor} +
\sum_{\lfloor n/2 \rfloor+1}^{n-1} \right]
\frac{n^{3/2}}{L(n)}\, q^{m*}(i) \, q(n-i) \, = \,
\\ \\ & & = \, \sum_{i=1}^{\lfloor n/2 \rfloor} q^{m*}(i) \, \bigg( \frac{n^{3/2}}{L(n)} \, q(n-i) \bigg) \,
+\, \sum_{i=1}^{\lceil n/2 \rceil - 1}
\bigg( \frac{n^{3/2}}{L(n)} \, q^{m*}(n-i) \bigg) \, q(i),
\end{eqnarray*}
and by dominated convergence the claim follows.

By the uniform convergence property of slowly varying sequences \cite[Th.~1.2.1]{cf:BinGolTeu}
we have that $L(c t)/L(t) \to 1$ as $t\to\infty$ {\sl uniformly in $c \in [\ga, 1/\ga]$}, for every
$\ga > 0$. In particular, we can choose $n_0$ such that
\begin{equation} \label{eq:unif_conv}
    L(ct) \le \sqrt 2 L(t) \qquad \quad \forall c \in [1/2, 1], \quad \forall t \ge n_0\,.
\end{equation}
Then we prove that there exists a constant $C >0$ such that:
\begin{equation}\label{eq2}
q^{k *}(n) \, \leq \, C \,k^3 \, \frac{L(n)}{n^{3/2}}, \qquad \forall k,n\in{\mathbb N},\ n \ge n_0.
\end{equation}
Again, we argue by induction: by~\eqref{c_q} the case $k=1$ holds for some positive
constant~$C$. If~\eqref{eq2} holds for all $k<2m$, $m\in\bbN$, then for $n\ge n_0$ we get:
\begin{eqnarray*}
q^{2m*}(n) & = & \sum_{i=1}^{\lfloor n/2 \rfloor}
q^{m *}(i) \, q^{m *}(n-i) \, + \sum_{i=\lfloor n/2 \rfloor +1}^{n-1}
q^{m *}(i) \, q^{m *}(n-i)
\\ \\ & \leq & 2 \sum_{i=1}^{\lfloor n/2 \rfloor} q^{m *}(i) \, q^{m *}(n-i)
\, \leq \, 2 \sum_{i=1}^{\lfloor n/2 \rfloor} q^{m *}(i) \,
C \, m^3 \, \frac{L(n-i)}{(n-i)^{3/2}}
\\ \\ & \leq &  2 \,C \, m^3\, \frac{\sqrt 2\, L(n)}{(n/2)^{3/2}} \sum_{i=1}^{\lfloor n/2 \rfloor} q^{m *}(i)
\, \leq \, C \, (2m)^3 \,\frac{L(n)}{n^{3/2}} \,.
\end{eqnarray*}
The case $k=2m+1$ follows similarly and consequently \eqref{eq2} holds. Therefore
we can apply dominated convergence in~\eqref{eq:rel1} (we recall that $\delta\in(0,1)$),
getting
\[
\frac{n^{3/2}}{L(n)} \, \widetilde Z^\rc_{\delta,n} \, = \,
\sum_{k=1}^\infty \, \frac{n^{3/2}}{L(n)} \, q^{k*}(n) \, \delta^k
\, \xrightarrow{n\to\infty} \, \sum_{k=1}^\infty k \, c_q \,  \delta^k \, = \,
\frac \delta{(1-\delta)^2} \ c_q.
\]
and \eqref{e3} is proven. Finally equation~\eqref{e32} follows by \eqref{rf} and \eqref{c_Y}
applying again dominated convergence:
\[
\frac{n^{1/2}}{L'(n)} \, \tilde Z_{\gd, n}^\rf \, = \, \sum_{t=0}^{n}
\tilde Z_{\gd, t}^\rc \ \frac{n^{1/2}}{L'(n)} \, P(n-t) \, \xrightarrow{n\to\infty} \,
\sum_{t=0}^\infty
\tilde Z_{\gd, t}^\rc \, = \, \sum_{k=0}^\infty \delta^k \, = \,
\frac{1}{1-\delta}\,,
\]
where in the second last equality we have used~\eqref{eq:ren_sol}.\qed

\subsection{The critical case.}
We treat now the case $\delta=1$.
\begin{lemma}\label{lemma4}
If $\delta=1$ and relations \eqref{c_q} and \eqref{c_Y} hold, then
\begin{gather}\label{e43}
\widetilde Z^\rc_{1,N} \, \sim \, \frac 1{2\pi} \, \frac{1}{c_q \, L(N) \, \sqrt{N}}
\qquad \quad (N\to\infty) \,,\\
\label{e432}
\widetilde Z^\rf_{1,N} \, \sim \, \frac{ L'(N)}{2 \, c_q \, L(N)} \qquad
\quad (N\to\infty)\,.
\end{gather}
\end{lemma}
\noindent
{\bf Proof}.
When $\gd = 1$ it is clear from~\eqref{eq:ren_sol} that $\widetilde Z^\rc_{1,N}$
is the Green function of the renewal process with step distribution~$q(\cdot)$.
More explicitly, if we set
\begin{equation}\label{renew}
\xi_k \, := \, T_1+\cdots+T_k\,, \qquad (T_i)_i \ {\rm IID}\,, \qquad
\p(T_i=n)=q(n)\,, \ \ n\in\N\,,
\end{equation}
then the law of $\xi_k$ is $q^{k*}(\cdot)$ and it is immediate to check that
$\widetilde Z^\rc_{1,N} = \bbP (\exists k:\ \xi_k = N)$. Then the asymptotic
behavior~\eqref{e43} is a result of Doney's \cite[Th.~B]{cf:Don97}.

To prove \eqref{e432}, we split the sum in~\eqref{rf} in three parts:
\begin{equation} \label{eq:split}
    \frac{L(n)}{L'(n)} \,\tilde Z_{1, n}^\rf \, = \, \Bigg(
    \sum_{t=0}^{\lfloor \ga n \rfloor} \ + \ \sum_{t= \lfloor \ga n \rfloor + 1}^{\lfloor (1-\ga) n\rfloor - 1}
    \ + \ \sum_{t=\lfloor(1- \ga) n \rfloor}^{n} \Bigg) \ \bigg( \,
    \frac{L(n)}{L'(n)} \, \widetilde Z_{\gd, t}^\rc \, P(n-t) \, \bigg)\,.
\end{equation}
Combining the asymptotic relations \eqref{c_Y} and \eqref{e43} with the uniform convergence
property of slowly varying sequences mentioned before equation \eqref{eq:unif_conv},
it is easy to check that the second sum above converges as~$n\to\infty$ to the integral
\begin{equation*}
    \frac{1}{2\pi \, c_q} \, \int_\ga^{1-\ga} \frac{dy}{\sqrt{y(1-y)}}\,.
\end{equation*}
On the other hand the fact that $\sum_{i=1}^k L(i)/\sqrt{i} \sim 2 L(k) \sqrt k$
as $k\to\infty$, cf. \cite[Prop.~1.5.8]{cf:BinGolTeu}, entails that the limits as $n\to\infty$
of the first and third sums in~\eqref{eq:split} are vanishing as~$\ga \to 0$, and
equation~\eqref{e432} follows.\qed

\subsection{The localized case.}
\label{S}
Let $\delta>1$. By continuity, there exists
$\tf_\delta>0$ such that
\[
\gd \sum_{t=1}^\infty
q(t) \, \exp (-\tf_\delta \, t) \, = \, 1.
\]
We set $q_\delta (t):= \delta \, q(t) \, \exp (-\tf_\delta t)$, $t\in\N$, so that $q_\delta(\cdot)$ is
a probability measure on $\N$. Notice that $\mu_\delta := \sum_t t \, q_\delta (t)<\infty$.
\begin{lemma}\label{lemma5}
If $\delta>1$ then
\begin{gather}\label{e5}
\widetilde Z^\rc_{\delta,N} \, \sim \, \frac 1{\mu_\delta} \, \exp(N \, \tf_\delta)
\qquad \quad (N \to\infty)\,,\\
\label{e52}
\widetilde Z^\rf_{\delta,N} \, \sim \, \Bigg( \frac 1{\mu_\delta} \,
\sum_{t=0}^\infty  e^{-\tf_\delta t} \, P(t) \Bigg)  \exp (N \, \tf_\delta)
\qquad \quad (N \to\infty)\,.
\end{gather}
\end{lemma}
\noindent
{\bf Proof}.
From~\eqref{eq:ren_sol} it is immediately seen that
\begin{equation*}
    e^{-\tf_\gd N} \, \widetilde Z^\rc_{1,N} \,=\, \sum_{k=0}^\infty
    q_{\gd}{}^{k*}(N)\,.
\end{equation*}
Arguing as in the proof of Lemma~\ref{lemma4} we have that the r.h.s. above
is the Green function of the renewal process with step distribution
$q_\gd(\cdot)$. However this distribution has \textsl{finite mean} $\mu_\gd$ and
therefore equation~\eqref{e5} is nothing but the standard Renewal Theorem~\cite{asmu}.

Finally, to prove~\eqref{e52} we resort to~(\ref{rf}):
\[
\exp(-\tf_\delta \, N) \, \tilde Z_{\gd, N}^\rf \, = \,
\sum_{t=0}^{N} \Big( e^{-\tf_\gd (N-t)} \tilde Z_{\gd, N-t}^\rc \Big)
e^{-\tf_\delta t} P(t) \, \xrightarrow{N\to\infty} \, \frac 1{\mu_\delta} \,
\sum_{t=0}^\infty e^{-\tf_\delta t} \, P(t)\,,
\]
having applied~\eqref{e5} and dominated convergence.\qed

\section{Infinite volume measures}
\label{sec:inf_vol}

Now we apply the asymptotic results obtained in the preceding section to the
{\sl thermodynamic limit} issue. Although $\bP_{\gep , N}^\rf$ and $\bP_{\gep , N}^\rc$
have been defined as measures on $(\R^+)^N$, it is convenient to extend them
to $(\R^+)^\N$ in an arbitrary way (for example by multiplying them by
$\prod_{i=N+1}^\infty \gd_0 ( \dd x_i)\,$). Our main result is the following theorem.

\begin{theorem}
\label{th:main}
For every $\gep\geq 0$ both $\bP_{\gep , N}^\rf$ and $\bP_{\gep , N}^\rc$ as measures on
$(\R^+)^\N$ converge weakly as $N \to \infty$ to the same limit $\bP_\gep$, law of
an irreducible Markov chain which is:
\begin{enumerate}
\item positive recurrent if $\gep  >\gep _c$ (localized regime)
\item transient if $\gep <\gep _c$ (strictly delocalized regime)
\item null recurrent if $\gep =\gep _c$ (critical regime)
\end{enumerate}
\end{theorem}

Let us introduce the times $(\tau_k)_{k\ge 0}$ at which the interfaces touches the wall:
\[
\tau_0:=0 \quad \qquad \tau_j :=\inf\{ n>\tau_{j-1}:\, x_n=0 \}
\qquad \quad x\in(\bbR^+)^\N\,,
\]
and the excursions $(e_k(\cdot))_{k \ge 0}$ of the interface above the wall:
\[
e_k(i) \, := \, \{ x_{\tau_k + i}:\, i=0,\ldots,\tau_{k+1}-\tau_k \}
\qquad \quad x\in(\bbR^+)^\N \,.
\]
We also set $\iota_N:=\sup\{k:\, \tau_k\le N\}$. The law of $(\tau_k)_{k\le\iota_N}$
under~$\bP_{\gep , N}^a$ can be viewed as a probability measure $p_{\gep , N}^a$
on the class $\cA_N$ of subsets of $\{1,\ldots, N\}$: indeed for $A\in\cA_N$, writing
\begin{equation}\label{notA}
A  = \{t_1,\ldots,t_{|A|}\}, \qquad
0 \, =: \, t_0<t_1<\cdots<t_{|A|} \, \leq \, N,
\end{equation}
we can set
\begin{equation}\label{defp^a}
p_{\gep , N}^a(A) \, := \, \bP_{\gep , N}^a(\tau_i= t_i, \
i\leq\iota_N).
\end{equation}
From the inclusion of $\cA_N$ into $\{0,1\}^\N$, the family of all subsets of~$\N$,
$p_{\gep , N}^a$ can be viewed as a measure on~$\{0,1\}^\N$.

The fundamental observation is that under $\bP_{\gep , N}^a$,
{\sl conditionally on $(\tau_k)_{k\le\iota_N}$}, the excursions $(e_k)_{0\leq k\leq \iota_N-1}$
are {\sl independent} and their laws are the same as under the unperturbed measure $\bP$.
Therefore, if we can prove that $p_{\gep , N}^a$ has a weak limit on $\{0,1\}^\N$,
then it is easily checked that the measure $\bP_{\gep, N}^a$ on $(\R^+)^\N$ converges weakly too,
and the limiting measure $\bP_\gep$ can be constructed simply by pasting the excursion over the
limit zero set. Therefore the analysis is completed by the following result.
\begin{prop}\label{str}
For every $\gep \ge 0$ both $p_{\gep , N}^\rf$ and $p_{\gep , N}^\rc$ as measures
on $\{0,1\}^\N$ converge to the same limit $p_\gep$, under which
$(\tau_k)_{k\geq 0}$ is:
\begin{enumerate}
\item for $\gep>\gep_c$ a renewal process with
interarrival probability $q_\delta(\cdot)$ (see \S \ref{S}).
\item for $\gep<\gep_c$ a terminating renewal process with
defective interarrival probability $\delta \, q(\cdot)$
\item for $\gep=\gep_c$ a renewal process with interarrival
probability $q(\cdot)$.
\end{enumerate}
\end{prop}

\medskip\noindent
{\bf Proof.} For all $\epsilon\geq 0$ and $a=\rf,\rc$, we have by the Markov property:
\[
\bP^a_{\gep, N} \left( \tau_1 =k_1, \tau_2=k_2,\ldots,
\tau_j=k_j \right)\, =\, \left[\prod_{i=1}^j \delta \, q (k_i-k_{i-1})\right]
\frac{\widetilde Z^a_{\delta,\, N+ \ind_{(a = \rc)} -k_j}}
{\widetilde Z^a_{\delta,\, N+ \ind_{(a = \rc)} }},
\]
for all $0=:k_0<k_1<\cdots<k_j\leq N$ (the factor $\ind_{(a = \rc)}$ in the partition functions
is due to definition~\eqref{eq:zeta_tilde}). Letting $N\to\infty$, we obtain
the thesis by Lemmas \ref{lemma2}-\ref{lemma4}-\ref{lemma5}. \qed

\subsection{More on the strictly delocalized regime}
If $\gep<\gep_c$, i.e. $\delta<1$, then under $\bP_\gep$ the number of
returns to $0$ and the last return to $0$ are a.s. finite random variables.
Their distributions are given in the following proposition, whose proof is
a straightforward consequence of Proposition~\ref{str}.
\begin{prop}\label{comp}
Let $\cN:=\#\{i\in\N: x_i=0\}$ and $\cL:=\sup\{i\in\N:
x_i=0\}$. Then for $\gep<\gep_c$:
\begin{equation}\label{comp1}
\bP_\gep(\cN = k) \, = \, (1-\delta) \, \delta^k, \qquad
k=0,1,\ldots
\end{equation}
\begin{equation}\label{comp2}
\bP_\gep(\cL = k) \, = \, (1-\delta) \, \tilde
Z^\rc_{\delta,k}, \qquad k=0,1,\ldots
\end{equation}
\end{prop}

\section{Scaling limits}
\label{sec:scaling_lim}

We finally turn to the scaling limits of our model.
We denote by $\big(X^N_t\big)_{t \in [0,1]}$ the linear interpolation of $\big(S_{i/N} \cdot L(N) /
\sqrt{N}\big)_{i=0,\ldots,N}$, the choice of the norming sequence being the natural one, see
Appendix~\ref{app:1}. We are interested in the weak convergence in $C([0,1])$
of the law of $\big(X^N_t\big)_{t \in [0,1]}$ under $\bP_{\gep, N}^a$.
This problem has been solved in~\cite{dgz} in the finite variance continuous set--up,
but for $\gep \neq \gep_c$ the techniques can be adapted in a straightforward way
to treat the general continuous and discrete settings considered here.

Consequently we focus on the critical case $\gep = \gep_c$. In fact in this regime
the result proven in~\cite{dgz} for $a=\rc$ is not optimal, the reason being
that the authors were not aware of Doney's result \cite[Th.~B]{cf:Don97}
which yields \eqref{e43}. In this section we show that the sharp asymptotic relations
\eqref{e43} and \eqref{e432} allow to simplify significantly the arguments in \cite{dgz},
proving the following

\begin{theorem}\label{ff} If $\delta=1$ then the process
$\big(X^N_t\big)_{t \in [0,1]}$ under $\bP_{\gep,N}^a$ converges in distribution
to the reflecting Brownian motion on $[0,1]$ for $a=\rf$ and to the reflecting
Brownian bridge on $[0,1]$ for $a=\rc$.
\end{theorem}

In the preceding section we have shown that, under the measure $\bP_{\gep, N}^a$,
there is a remarkable decoupling between the {\sl zero level set} (that is the set of
points where the interface touches the wall) and the excursions of the interface above
the wall. Namely, conditionally on the zero level set, the excursions are an independent
family and their laws are the same as under the initial measure $\bP$. Since our basic assumption
(H) entails that the measure $\bP$ is attracted to the Gaussian law, it is not a surprise
that the law of the rescaled excursion under $\bP$ converges weakly to the law of the {\sl Brownian
excursion}. The proof of this fact can be found for instance in~\cite{dgz} for the continuous
set--up and in \cite{cf:K} for the discrete set--up (these proofs are given in the case
of finite variance, but they can be easily adapted to the general case).

Therefore we focus our attention to the law of the rescaled zero level
set. More precisely we introduce
${\cA}_N^a$, a random subset of $[0,1]$, by setting $\mathbb P \left( {\cA}^a_N =A/N\right)=p_{\gep _c
,N}^a(A)$, $a=\rf,\rc$, for $A \subseteq \{0, \ldots, N \}$ (recall the notation introduced in
\eqref{defp^a}). If we can prove the weak convergence of the zero set ${\cA}_N^a$,
then, in view of the convergence of the excursions mentioned above, the weak
convergence of the full measure $\bP_{\gep,N}^a$ follows arguing like in Section 8 of~\cite{dgz}.
Therefore Theorem~\ref{ff} is a consequence of the following proposition, first proven
in \cite[Prop. 10]{dgz} (in the finite variance continuous case and
with an additional assumption for $a=\rc$).
\begin{prop}\label{c0bm}
As $N\to \infty$ we have that:
\begin{itemize}
\item[(i)]
${\cA}_N^\rf$ converges in law to
$\left\{t\in [0,1]: B(t)=0 \right\}$,
\item[(ii)]
${\cA}_N^\rc$ converges in law to
$\left\{t\in [0,1]: \beta(t)=0 \right\}$,
\end{itemize}
where $B$ is a standard Brownian motion and $\beta$ is a Brownian bridge over $[0,1]$.
\end{prop}
\noindent
The basic notions about the convergence in law of random sets are recalled in Appendix~\ref{sec:j-d}
(for more details see \cite[\S~3]{cf:FFM} and~\cite{cf:Matheron}).

It is convenient to introduce a simpler random set $\cA_N$,
to which the random sets ${\cA}_N^\rf$ and ${\cA}_N^\rc$ are strictly linked.
Namely we consider again the renewal process $\xi_k=T_1+\cdots+T_k$, where $(T_i)_i$ is IID
and $\p(T=n)=q(n)$, $n\in\N$, and we set ${\cA}_N:=\{\xi_k/N: k\in\N\} \cap [0,1]$.
Then the asymptotic relation \eqref{c_q} for $q(\cdot)$ implies
the following basic result, first proven in~\cite[Lemma 5]{dgz}
for the case in which $L(\cdot)$ is a constant (but the proof
extends to our general set--up in a straightforward way).
\begin{lemma}\label{th:JB}
The sequence $(\cA_N)_N$ converges in law to $\{ t \in [0,1]: B(t) =0\}$.
\end{lemma}
\noindent
We point out that the proof given in~\cite{dgz} uses in an essential way the theory
of {\sl regenerative sets} and their connection with {\sl subordinators} (we refer to~\cite{cf:FFM}
for more on this subject). In view of the importance of this result, in
Appendix~\ref{sec:j-d} we sketch an alternative and more direct proof,
which is built on the basic relation~\eqref{e43}.

\medskip
\noindent
{\bf Proof of Proposition~\ref{c0bm}.}
Let us first consider the free case $a=\rf$.
It is easy to see that the laws of ${\cA}_N^\rf$ and ${\cA}_N$ are equivalent, more
precisely for every bounded measurable functional $\Phi$ we have
\[
\bbE\Big[\Phi({\cA}_N^\rf)\Big] =\,
\bbE\Big[\Phi({\cA}_N) \, f^\rf_N(\sup {\cA}_N)\Big]\,, \qquad
f^\rf_N(t) \, := \,  \frac{ P \left(N \, (1-t)\right)}
{{\widetilde Z^\rf_{1,N}} \cdot Q \left( N \, (1-t)\right)}\,, \quad t\in[0,1]\,,
\]
where $Q(n):=\sum_{t=n+1}^\infty q(t)$.
The asymptotic behavior of $q(\cdot)$ being given by~\eqref{c_q}, it follows from
\cite[Prop.~1.5.10]{cf:BinGolTeu} that $Q(n) \sim 2 c_q L(n) / \sqrt{n}$ as $n\to\infty$.
Hence by \eqref{c_Y} and \eqref{e432} one sees that $\lim_{N\to\infty} f^\rf_N(t) = 1$
uniformly in $t\in[0,\gamma]$, for every $\gamma\in(0,1)$, and then (i)
is an easy consequence of Lemma~\ref{th:JB}.

We turn now to $a=\rc$. Here it is more convenient to study the Radon-Nikodym derivative
of the law of ${\cA}_N^\rc\cap [0,1/2]$ w.r.t. the law of
${\cA}_N\cap [0,1/2]$. This time the Radon--Nikodym derivative is given
by (cf. \cite[Proof of (32), Step 1]{dgz})
\[
\bbE\left[\Phi({\cA}_N^\rc\cap [0,1/2])\right] \, = \,
\bbE\left[\Phi({\cA}_N\cap [0,1/2]) \ f^\rc_N\Big(\sup \big( {\cA}_N\cap [0,1/2] \big) \Big)\right]\,,
\]
\[
f^\rc_N(t) \, := \,  \frac{\sum_{n=0}^{N/2} {\tilde Z^\rc_{1,n}} \, q(N+1-Nt-n)}
{{\tilde Z^\rc_{1,N+1}} \ Q(\lfloor N/2 \rfloor-Nt)}\,, \quad t\in[0,1/2]\,.
\]
By \eqref{c_q}, \eqref{c_Y} and \eqref{e43} we see that:
\[
\lim_{N\to\infty} f^\rc_N(t) \, = \,  \frac {\int_0^{1/2} y^{-1/2}\, (1-t-y)^{-3/2} \, dy}
{2 \, (1/2-t)^{-1/2}} = \frac{1}{\sqrt{2}} \, \frac{1}{1-t}\,,
\]
uniformly in $t\in[0,\gamma]$, for any $\gamma\in(0,1/2)$. Then (ii) follows
from Lemma~\ref{th:JB} by the same arguments used in \cite[Proof of (32), Step 3]{dgz}. \qed

\appendix

\section{An asymptotic relation}
\label{app:asympt}

We are going to prove that relation \eqref{eqq} holds true.

\subsection{A Local Limit Theorem}
\label{app:1}
We recall that, by our basic assumption~(H), one has the weak convergence
\begin{equation} \label{eq:clt}
    \text{under $\bP\,$:} \qquad \frac{L(n)}{\sqrt n} \, S_n \, \Rightarrow \, \frac{1}{\sqrt{2\pi}}
    \, e^{-x^2/2} \, \dd x \qquad \quad (n\to\infty)\,,
\end{equation}
where $L(\cdot)$ is a slowly varying function satisfying the
relation $L(x) \sim 1/ \sqrt{\cV(\sqrt{x} / L(x))}$ as~$x \to \infty$. More explicitly,
this function can be defined as $L(x) := \sqrt{x} / g^{-1}(x)\,$,
where $g(\cdot)$ is any {\sl increasing}
function such that $g(x) \sim x^2 / \cV(x)$ as~$n\to\infty$ (the existence of such $g(\cdot)$
is guaranteed by~\cite[Th.~1.5.3]{cf:BinGolTeu}, where an explicit definition is given).

We point out that equation~\eqref{eq:clt} expresses the most general instance in which a
random walk is attracted (without centering) to the Gaussian law, which in turn happens
if and only if condition (H) holds, cf.~\cite[\S IX.8 \& \S XVII.5]{cf:Feller2}.
Of course in the special case $\sigma ^2 := \bE\big[|Y_1|^2\big]<\infty$ we have
$L(t) \equiv \gs^{-1}$ by the Central Limit Theorem.

\smallskip

Let us denote by $f_n(x)$ the density (resp. the mass function) of $S_n$ under~$\bP$,
in the continuous set--up (resp. in the discrete set--up). Then the Local Limit Theorem for Densities,
cf.~\cite[\S 46]{cf:GneKol}, (resp. Gnedenko's Local Limit Theorem,
cf.~\cite[\S 8.4]{cf:BinGolTeu}) yields the asymptotic relation
\begin{equation} \label{eq:llt}
    f_n(0) \,\sim\, \frac{1}{\sqrt{2\pi}} \, \frac{L(n)}{\sqrt n} \qquad \quad (n\to\infty)\,.
\end{equation}

\subsection{The continuous case} We follow the proof given in~\cite{dgz}
in the case of finite variance. Introducing the set
$\cC_n := \{x_1 > 0, \ldots, x_n > 0\}$, the very definition of $Z_{0, n}^\rc$ gives
\begin{equation} \label{eq:app1}
    Z_{0,n}^\rc \,=\, \int_{\cC_n} e^{-H_n^\rc(x_1, \ldots, x_n)} \,
    \dd x_1 \cdots \dd x_n \,.
\end{equation}
On the other hand for the density of $S_{n+1}$ under $\bP$ we have
\begin{equation*}
    f_{n+1}(0) \,=\, \int_{\R^n} \frac{e^{-H_n^\rc(x_1, \ldots, x_n)}}{\kappa^{n+1}} \,
    \dd x_1 \cdots \dd x_n \,.
\end{equation*}
Next we introduce the linear transformation $T_n: \R^n \to \R^n$ defined by
\begin{equation*}
    T_n(x_1, \ldots, x_n) \; := \; (x_2 - x_1,\, x_3 - x_1,\, \ldots,\, x_n - x_1,\, -x_1)\,,
\end{equation*}
which is nothing but a cyclical permutation of the increments of $(0, x_1, \ldots, x_n ,0)$.
Notice that $T_n$ preserves Lebesgue measure, that $(T_n)^{n+1}$ is the identity map on $\R^n$
and that $H_n^\rc$ is invariant along the orbits of $T_n$, namely $H_n^\rc(T_n (x)) = H_n^\rc(x)$.
Moreover, the $n+1$ sets $\{(T_n)^k(\cC_n)\,, \ k=0, \ldots, n\}$ are {\sl disjoint} and their
union differs from the whole $\R^n$ only by a set of zero Lebesgue measure.
These considerations yield
\begin{equation*}
    f_{n+1}(0) \,=\, \sum_{k=0}^n \, \int_{(T_n)^k(\cC_n)} \frac{e^{-H_n^\rc(x)}}{\kappa^{n+1}}
    \, \dd x \,=\, (n+1) \int_{\cC_n} \frac{e^{-H_n^\rc(x)}}{\kappa^{n+1}} \,\dd x\,,
\end{equation*}
and comparing with \eqref{eq:app1} we get
\begin{equation*}
    \frac{Z_{0,n}^\rc}{\kappa^{n+1}} \,=\, \frac{1}{n+1} \, f_{n+1}(0)\,.
\end{equation*}
Therefore it suffices to apply~\eqref{eq:llt} to show that relation
\eqref{eqq} holds with~$C=1$.

\subsection{The discrete case} \label{app:discrete}
For $S=(S_n)_{n\ge 0} \in \R^\N$ we introduce the
{\sl weak descending} ladder epochs
\begin{equation*}
    \overline T_0 := 0 \qquad \quad \overline T_{k+1} := \inf \{ n > T_k:\ S_n \le S_{T_k} \}
\end{equation*}
and the corresponding ladder heights $\overline H_k := - S_{T_k}$. The reason why these quantities are of
interest to us is that the definition of $Z_{0, n}^\rc$ yields
\begin{equation} \label{eq:app2}
    \frac{Z_{0,n}^\rc}{\kappa ^{n+1}} \,=\, \bP \big( \overline T_1 = n+1,\,\overline H_1 = 0 \big) \,.
\end{equation}
Now we are going to use a fundamental combinatorial identity discovered by Alili and Doney,
that for $k,\, n \in \N$ and $x \ge 0$ reads as
\begin{equation*}
    \bP \big( \overline T_k = n,\,\overline H_k = x \big) \,=\, \frac{k}{n} \,
    \bP \big( \overline H_{k-1} \le x < \overline H_k,\, S_n = x \big)\,,
\end{equation*}
cf. \cite[Eq.~(3)]{cf:AliDon} (the interchange between $<$ and $\le$ with respect to that
formula is due to the fact that they consider {\sl strong} instead of weak ladder variables).
Plugging this identity into \eqref{eq:app2} and arguing as in the step~(i) of the proof of
Proposition~6 in~\cite{cf:AliDon}, we get
\begin{equation*}
    \frac{Z_{0,n}^\rc}{\kappa ^{n+1}} \,=\, \frac{1}{n+1} \,
    \bP \big( \overline H_1 > 0,\, S_n = 0 \big) \,\sim\,
    \frac{1}{n} \, \bP \big( \overline H_1 > 0 \big) \, \bP \big( S_n = 0 \big)
    \,\sim\, \frac{\bP \big( \overline H_1 > 0 \big)}{\sqrt{2\pi}} \, \frac{L(n)}{n^{3/2}}\,,
\end{equation*}
where we have applied \eqref{eq:llt}, and equation \eqref{eqq} is proved with $C=\bP(\overline H_1 > 0)$.

\section{The critical zero level set}
\label{sec:j-d}

We want to sketch here an alternative proof of Lemma~\ref{th:JB}. For the purpose of
this section it is convenient to consider the random set $\cA_N$ of Section~\ref{sec:scaling_lim}
on the whole positive real line instead of the interval $[0,1]$. More precisely, introducing
the renewal process $\xi_k=T_1+\cdots+T_k$, where $(T_i)_{i \in \N}$ is IID and $\bbP(T_1=n)=q(n)$,
we set ${\cA}_N:=\{\xi_k/N: k\in\N\}$.

\smallskip

Let us first recall
some basic facts on the convergence of closed sets. We denote by $\cF$ the family of all
closed sets of $\R^+$, and we endow it with the topology of Matheron,
cf.~\cite{cf:Matheron} and \cite[\S~3]{cf:FFM}, which in our setting can be conveniently
described as follows. For $F \in \cF$ and $t\in\R^+$ we set $d_t(F) := \inf\big(F \cap (t,\infty)\big)$.
Notice that $t \mapsto d_t(F)$ is a right--continuous function and that
the set $F$ can be actually identified with the function
$d_{(\cdot)}(F)$, because $F = \{t\in\R^+:\,d_{t-}(F) = t\}$.
Then in terms of $d_{(\cdot)}(F)$ the Matheron topology is the standard Skorohod topology
on c\`adl\`ag functions taking values in
$\overline{\R^+} := \R^+ \cup \{+\infty\}$. We point out that with this topology the space $\cF$ is metrizable,
separable and {\sl compact}, hence in particular Polish. Moreover the Borel $\gs$--field
on $\cF$ coincides with the $\gs$--field generated by the maps $\{d_t(\cdot),\, t\in\R^+\}$.

\smallskip

Let us denote respectively by $\bbP_N$ and $\bbP^{(BM)}$ the laws of the random closed sets
$\cA_N$ and $\{t\in\R^+:\, B(t) = 0\}$, where $B(\cdot)$ is a
standard Brownian motion. These laws are probability measure on $\cF$, and our goal is to prove
that $\bbP_N$ converges weakly to $\bbP^{(BM)}$ as $N\to\infty$.
Thanks to the compactness of $\cF$, we can focus on
the convergence of the {\sl marginal distributions}. More precisely,
it is sufficient to show that for every $n\in\N$ and for all $t_1 < \ldots < t_n \in \R^+$
one has the weak convergence of the image laws on $(\R^+)^n$:
\begin{equation} \label{eq:weak}
    \bbP_{N} \circ (d_{t_1}, \ldots, d_{t_n})^{-1} \ \Longrightarrow \
    \bbP^{(BM)} \circ (d_{t_1}, \ldots, d_{t_n})^{-1} \qquad (N\to\infty)\,,
\end{equation}
and the result $\bbP_N \Rightarrow \bbP^{(BM)}$ follows, because the
distributions of the vectors $(d_{t_1}, \ldots, d_{t_n})$ determine laws on $\cF$.

The validity of \eqref{eq:weak} can be obtained
by direct computation. For simplicity we will only consider the case $n=1$,
the general case follows along the same line.
We recall that for any~$t>0$, the law of~$d_t$ under~$\bbP^{(BM)}$ is given by
\begin{equation*}
    \bbP^{(BM)} \big( d_t \in \dd y \big) \;=\; \frac{t^{1/2}}{\pi \, y(y-t)^{1/2}} \,
\ind_{(y > t)} \, \dd y \;=:\; \rho_t(y) \, \dd y \,,
\end{equation*}
cf.~\cite{reyo}, hence we have to show that for every~$x\in\R^+$
\[
\lim_{N\to\infty} \, \bbP_N \big(d_t\ge x\big)
\, = \,\int_x^\infty \rho_t(y) \, dy \,.
\]
Using the Markov property for the renewal process $(\xi_k)$ we get
\begin{align*}
& \bbP_N \big(d_t \ge x\big) = \sum_{k\in\N} \; \bbP \big(\xi_k \leq Nt\,,\;
\xi_{k+1}\ge Nx\big) \\
&\qquad = \sum_{i=1}^{\lfloor Nt \rfloor } \; \sum_{j= \lceil Nx \rceil}^\infty  \Bigg( \sum_{k\in\N} \;
\bbP \big( \xi_k=i \big) \Bigg)\, \bbP \big( \xi_1 = j-i \big) \\
&\qquad = \sum_{i=1}^{\lfloor Nt \rfloor }  \Bigg( \sum_{k\in\N} \;
q^{* k} ( i ) \Bigg) \sum_{j= \lceil Nx \rceil}^\infty q(j-i)
= \sum_{i=1}^{\lfloor Nt \rfloor } \; \widetilde Z_{1,i}^\rc
\; Q\big( \lceil Nx \rceil - i -1 \big)\,,
\end{align*}
where we have applied \eqref{eq:ren_sol}. We recall the notation $Q(n) := \sum_{k=n+1}^\infty q(k)$,
introduced in the proof of Proposition~\ref{c0bm}, and the fact that $Q(n) \sim 2 c_q L(n) / \sqrt{n}$
as $n\to\infty$, as it follows from \eqref{c_q} applying \cite[Prop.~1.5.10]{cf:BinGolTeu}.
Since the asymptotic behavior of the constrained partition function in the critical case is given
by~\eqref{e43}, we obtain
\begin{align*}
& \bbP_N \big(d_t \ge  x\big) \;\sim\; \sum_{i=1}^{\lfloor Nt \rfloor } \;
\frac{1}{2\pi} \, \frac{1}{c_q \, L(i) \, \sqrt{i}} \;
\frac{2 c_q L\big( \lceil Nx \rceil - i -1 \big)}{\sqrt{ \lceil Nx \rceil - i -1 }}
\qquad \ \ (N\to\infty)\,.
\end{align*}
Now using the fact that $L(ct)/L(t) \to 1$ as $t\to\infty$ {\sl uniformly} in $t \in [\ga, 1/\ga]$,
for every $\ga > 0$ (cf.~\cite[Th.~1.2.1]{cf:BinGolTeu}), and the convergence of the
Riemann sums to the corresponding integral we easily get
\begin{align*}
    & \exists \lim_{N\to\infty} \bbP_N \big(d_t \ge x\big) \;=\; \frac{1}{\pi}
    \int_0^{t} \dd s \, \frac{1}{\sqrt{s}} \, \frac{1}{\sqrt{x-s}}
    \;=\; \int_x^\infty \dd z \, \rho_t(z) \;,
\end{align*}
that is what was to be proven.

\end{document}